\documentclass[11pt]{article}
\usepackage{amsfonts} 
\usepackage{amssymb}
\usepackage{amsmath}  \usepackage{url} \usepackage{amsthm}
\usepackage{epsfig,epic,eepic}  \usepackage[usenames,dvips]{color}

\usepackage{psfrag}

\newcommand{\R}{\mathbb R}

\newcommand{\hop}{\vskip .2cm\noindent}
\newcommand{\hip}{\vskip .1cm\noindent}

\newcommand{\hD}{\widehat D}

\newcommand{\hM}{\widehat M}
\newcommand{\hg}{\widehat g}
\newcommand{\hX}{\widehat X}
\newcommand{\hY}{\widehat Y}\newcommand{\hT}{\widehat T}
\newcommand{\hZ}{\widehat Z}
\newcommand{\wX}{X}
\newcommand{\wY}{Y}
\newcommand{\wZ}{Z}
\newcommand{\wL}{{\mathcal{L}}}

\newcommand{\weg}[1]{}
\def\tr{{\textrm{\sl{trace}}}} 

\def\slsf#1{{\slshape \sffamily #1\/}}
\newtheorem{enonce}{}[section]

\newtheorem{thm}{Theorem}
\newtheorem{cor}[enonce]{Corollary}
\newtheorem{prop}[enonce]{Proposition}
\newtheorem{rema}[enonce]{Remark}
\newtheorem{fact}[enonce]{Fact}

\def\const{{\sf{const}}}

\author {Vladimir S. Matveev\thanks{ partially supported by  DFG (SPP 1154 and GK 1523)} \ and \  Pierre Mounoud}
\title{ Gallot-Tanno Theorem for closed incomplete pseudo-Riemannian manifolds and applications }
\date{}
\begin{document}
\maketitle
\begin{abstract}
We extend the Gallot-Tanno Theorem to closed pseudo-Riemannian manifolds. It is done by showing that if the cone over a manifold  admits a parallel symmetric $(0,2)-$tensor then it is Riemannian. Applications of this result  to the existence of metrics with distinct Levi-Civita connections but having the same unparametrized geodesics and to the projective Obata conjecture  are given.  We also apply our result to show that the holonomy group of a closed  $(O(p+1,q),S^{p,q})$-manifold does not preserve any nondegenerate splitting of $\R^{p+1,q}$.
\end{abstract}

\section{Introduction.} \subsection{Main result} 
Let $(M,g)$ be a pseudo-Riemannian manifold and let $D$ be the corresponding   covariant derivative. Along this paper we will consider Riemannian metrics  as pseudo-Riemannian.
We are interested in the existence of non-constant functions $\alpha$ on $M$ such that for all   vector fields  $X$, $Y$, $Z$ on  $M$ we have:
\begin{equation}\label{one}\qquad DDD \alpha (X,Y,Z)+ c\cdot \bigl( 2 (D\alpha  \otimes g)(X,Y,Z)+(D\alpha \otimes g)(Y,X,Z)+(D\alpha \otimes g)(Z,X,Y)\bigr)=0,\end{equation}
where $c$ is a real constant.  This is a linear system of  PDE on $\alpha$, in  the ``tensor" notation it reads 
$$
\alpha_{,ijk} = c\cdot \left(2 \alpha_{,k} g_{ij}  + \alpha_{,j} g_{ki} + \alpha_{,i} g_{kj}\right).
$$

The main result of this article is the following 
\begin{thm}\label{pr}
 Let  $(M,g)$ be a closed (i.e., compact without boundary) connected pseudo-Riemannian  manifold and $c$ be a real number. If $\alpha:M\to \R$ is a non-constant function  satisfying \eqref{one}, then 
$c\ne 0$ and  the metric $c\cdot g$ is {a} Riemannian (i.e., positively definite) metric of constant 
curvature $1$.
\end{thm}

\subsection{History,  motivation,  and applications} \label{his} 
The equation  \eqref{one} has already been studied,  mostly in the Riemannian setting. 
The motivation of Gallot and Tanno to study this equation came from the  spectral geometry:  it is well-known  (see for example \cite{Ga})  that, on the standard sphere $S^n \subset \mathbb{R}^{n+1}$ of dimension $n>1$,  all 
eigenfunctions corresponding to the second  biggest eigenvalue $-n$   of the laplacian satisfy the equation 
\begin{equation} \label{obata1}
 DD\alpha(X,Y) +  \alpha \cdot  g(X,Y)=0, \end{equation} (for all vector fields  $X$, $Y$, $Z$ on  $M$). The  
eigenfunctions corresponding to the  third  biggest eigenvalue $-2(n+1)$ satisfy  \eqref{one} with $c=1$. 

Moreover, Obata has shown \cite[Theorem A]{obata1} that, on closed  Riemannian  manifolds, 
 the existence of a nonconstant solution of \eqref{obata1}   implies that the metric has constant curvature $1$.   Later, he \cite{Obata}, and, according to Gallot \cite{Ga},  Lichnerowicz, 
  asked the question whether the same holds for the equation \eqref{one} (assuming $c=1$). The affirmative answer was given  in \cite{Ga,Ta}.  Theorem \ref{pr}, which is the main result of our paper,  generalizes the results of    \cite{Ga,Ta} to pseudo-Riemannian metrics.

This equation also appears in the context of geodesic equivalence. Recall that two metrics $g$ and $ g'$ are \slsf{geodesically equivalent} (\slsf{affinely equivalent}, resp.) 
 if every $g$-geodesic is a reparametrized $g'$-geodesic (if their Levi-Civita connections coincide, resp.) A vector field $V$ on $(M, g)$  is called \slsf{projective} (\slsf{affine}, resp.), if its local flow acts by  \slsf{projective transformations} (\slsf{affine transformations}, resp.), i.e.,  takes unparameterized geodesics to geodesics (preserves  the Levi-Civita connection, resp.).  

   In particular,  Solodovnikov \cite{So} has  showed that if a Riemannian metrics $g$ admits  ``lots''  (we formalize this notion in section  \ref{lichne}) of geodesically equivalent, but not affinely equivalent  metrics, 
 then there exists a real number $c$ and a non-constant solution of $(\ref{one})$, see section \ref{lichne} for definitions and more precise statements. This result has been recently extended to the pseudo-Riemannian setting  in Kiosak et al \cite{Ma2}.  Hiramatu \cite{hiramatu} has shown that if a Riemannian metric of constant scalar curvature on a closed manifold  admits  a nonaffine  projective vector field, then   there also exists a nonconstant solution of the equation $(\ref{one})$ for a certain constant $c$.
    Kiosak et al \cite{einstein}  has shown that if an Einstein metric $g$ admits  a geodesically equivalent, but not affinely equivalent  metric $g'$, then there exists a  non-constant solution of $(\ref{one})$ for a certain constant $c$. 
 
 The equation \eqref{one} naturally appears in the study of the geometry of the metric cones, see Gallot \cite{Ga} or  Alekseevsky et al \cite{Leist}.  We will explain the relation  between cones with decomposable  holonomy and the equation \eqref{one} in section \ref{section}. In  fact, this  relation is one of the main tools of our proof.

Combining theorem \ref{one}  with the results listed above, 
we obtain  that on a closed manifold  any pseudo-Riemannian metric admitting  ``lots" of geodesically  equivalent, but not affinely equivalent metrics    is, up to multiplication by a constant,  the Riemannian metric of constant curvature $1$ (cf. corollary \ref{corobata}).

  We also obtain 
 that Einstein pseudo-Riemannian metrics of nonconstant curvature  on closed manifolds are geodesically rigid, in the sense that every metric geodesically equivalent to   them are actually affinely  equivalent to them (cf. corollary \ref{ein}). 
 
   We also obtain that metric cones with decomposable holonomy  group over closed pseudo-Riemannian metrics are Riemannian and  flat (cf. proposition \ref{deco}).  {We  apply this result to closed  nonzero constant curvature (non-Riemannian) manifolds i.e.,  to manifolds locally modeled on a pseudosphere $S^{p,q}$. We obtain that, for such a manifold,  the holonomy group of the associated $(O(p+1,q),S^{p,q})$-structure does not preserve any non degenerate splitting of $\R^{p+1,q}$ (cf. corollary \ref{ghz}). In \cite{ze}, Zeghib proved this statement under the additional condition of completeness.}

\subsection{Previous results} 
Partial  versions  of theorem \ref{pr} were known before. In the Riemannian case, theorem  \ref{pr} is due to Gallot \cite[Corollaire 3.3]{Ga} and Tanno \cite[Theorem A]{Ta} under the assumption  $c>0$, and is due to Hiramatu \cite[Lemma 2]{hiramatu}  under the assumption  $c\le 0$.  

Moreover, Gallot and Tanno assumed only completeness (instead of closedness). In the realm of Riemannian geometry, the unit tangent bundle of a compact manifold being compact, closeness implies  completeness. This is no more true in the pseudo-Riemannian geometry, where  incomplete metrics on compact manifolds are abundant. For example,  by    Carri\`ere et al  \cite{CR} the set of incomplete Lorentzian $2$ dimensional tori is dense in the set of Lorentzian tori.  Completeness and closeness are quite independent properties in the pseudo-Riemannian geometry, and it is not an easy task to understand whether a given metric on a closed manifold  is complete. 
Moreover,  \cite[Example 3.1]{Leist} from Alekseevsky et al provides   non-compact complete pseudo-Riemannian manifolds of non-constant curvature admitting non-constant solutions to $(\ref{one})$.
Moreover, under the additional  assumption that the metric is complete theorem \ref{pr} is easy, see  \cite[Theorems 1,2]{Ma}. 

\subsection{Organisation of the paper and the converse statement.}
The round sphere  $S^n := \{(x^1,...,x^{n+1}) \in \mathbb{R}^{n+1}\mid (x^1)^2+...+ (x^{n+1})^2=1\}$ with the standard metric  admits a lot of nonconstant solutions of \eqref{one} (with $c=1$): as we mentioned in \S \ref{his}, every eigenfunction of the Laplacian  corresponding to the  third biggest
eigenvalue $-2(n + 1)$ satisfies  (1). { By our theorem,  any closed manifold admitting a non constant solution of \eqref{one} is,   up to  a constant, a quotient of $S^n $,  but certain of those quotients do not admit nonconstant solutions of \eqref{one}.

Indeed, let  $M$ be the quotient of $S^n $ by   a discrete subgroup $\Gamma\subset O(n+1)$. The cone over $M$ is the quotient of the cone over $S^n$ (i.e., of $\R^{n+1}\setminus \{0\}$ endowed with the euclidean metric) by $\Gamma$. By proposition \ref{para->decomp},  $M$ admits a non-trivial solution of \eqref{one}, if and only if its cone is decomposable i.e., if and only if $\Gamma$ preserves an orthogonal splitting of $\R^{n+1}$, or equivalently if and only if there exists $0<p<n+1$ such that $\Gamma \subset O(p)\times O(n+1-p)$. 
Thus, the only quotients of the sphere $S^3$ admitting nontrivial solutions of \eqref{one} are the lens spaces. It follows that  the Poincar\'e homology sphere (which is the quotient of the standard $3$-sphere 
by the lift of the group of direct isometries of the regular dodecahedron)
admits no nonconstant solution of \eqref{one}.
}

The organization of the article is as follows. In section \ref{sec0} we prove theorem \ref{pr} under the additional assumption $c=0$.  The rest of the paper is devoted to the case $c\ne 0$ ---  we will explain  in remark \ref{tttt},  that if  $c\ne 0$, then without loss of generality we can assume $c=1$. 
 In section \ref{section}  we establish a  link between solutions of $(\ref{one})$ (with $c=1$) and parallel symmetric $(0,2)-$tensors on the cone over $(M,g)$: we  show that the existence of a non-constant solution of $(\ref{one})$ is equivalent to  that the cone is decomposable. In section \ref{decompo}  decomposable  cones   are  studied and theorem \ref{pr} is proved.  Section \ref{lichne} is devoted to the application of theorem \ref{pr} in the theory of geodesically equivalent metrics. {Section \ref{appendix} is devoted to the study of the holonomy of closed manifolds with constant nonzero curvature.}

\section{Proof of theorem \ref{pr} under the assumption $c=0$.}  \label{sec0}

Assume   $c=0$. Equation $(\ref{one})$ implies that the Hessian of  $\alpha$ is  parallel. Since the manifold is closed,  $\alpha$ has a minimum and a maximum.
At a minimum,  the Hessian must be nonnegatively  definite, 
 and at a maximum it must be nonpositevely definite. Therefore the Hessian is null, and the gradient of $\alpha$ is parallel. But as it vanishes at the extremal points, it vanishes everywhere and $\alpha$ is constant. Theorem \ref{pr} is proved under the assumption $c=0$. 
 
 \begin{rema}  \label{tttt} If $c\ne 0$,  without loss of generality we can assume  $c=1$. Indeed, if a  function   $\alpha$ is a solution of \eqref{one} with $c\ne 0$, then it is also a solution of the equation 
 \begin{equation}  DDD \alpha (X,Y,Z)+  \bigl( 2 (D\alpha  \otimes g')(X,Y,Z)+(D\alpha \otimes g')(Y,X,Z)+(D\alpha \otimes g')(Z,X,Y)\bigr)=0 \label{two}\end{equation} for $g':= c\cdot g$.
  Since the Levi-Civita connections of $g$ and of $g'$ coincide, the equation \eqref{two} is the equation \eqref{one}  with respect to  the metric $g'$  with $c=1$.    \end{rema}

\section{Parallel symmetric $(0,2)-$tensors on the cone over a manifold and nonconstant solutions of \eqref{one} for $c=1$.}\label{section}
Let $(M,g)$ be a pseudo-Riemannian manifold. The  \slsf{cone manifold} over $(M,g)$ is  the manifold $\hM=\R_{>0}\times M$ endowed with the metric $\hg$ defined by $\hg=\,dr^2+ r^2 g$ (i.e., in the local coordinate system $(r,x^1,...,x^n) $ on $\hM$, where $r$ is the standard coordinate on $\R_{>0}$, and $(x^1,...,x^n)$ is a local coordinate system on $M$, the scalar product in $\hg$ of the vectors 
$u= u_0 {\partial_r} + \sum_{i=1}^n u^i {\partial_{ x^i} } $ and $v= v_0 {\partial_r} + \sum_{i=1}^n v^i {\partial_{x^i} } $ is given by $\hg(u,v)= u^0v^0+ r^2\sum_{i=1}^n g_{ij} u^i v^i$).

We will denote by $D$ the Levi-Civita connection of $g$ and by $\hD$ the Levi-Civita connection of $\hg$.

 The holonomy of cones over pseudo-Riemannian is strongly related to the equation (\ref{one}). This relation is given by the following proposition, which is almost contained in the proofs of   \cite[corollaire 3.3]{Ga} (for an implication)  and   in \cite[Corollary 1]{Ma} (for the reciprocal). As we will use  some lines from it, as those proofs have a non empty intersection,  and for the convenience of the reader, we give its proof but it does not pretend to be new.
\begin{prop}\label{GaMa}
Let $(M,g)$ be a pseudo-Riemannian manifold, let $c=1$. 
Let  $(\hM,\hg)$ be the cone manifold over $(M,g)$.

Then, there exists a \slsf{non-constant}  function $\alpha: M\rightarrow \R$ satisfying \eqref{one},   
if and only if there exists a \slsf{non-trivial} (i.e.,  not proportional to $\hg$) symmetric parallel  (i.e., the covariant derivative vanishes) $(0,2)-$tensor on $(\hM,\hg)$.

More precisely if $\alpha$ is a non-constant  solution of $(\ref{one})$ then the Hessian of the function $A:\hM\rightarrow \R$ defined by $A(r,m)=r^2\alpha(m)$ is non-trivial and parallel (i.e., $\hD\hD\hD A=0$). Conversely if $\hT$ is a non-trivial symmetric parallel $(0,2)-$tensor on $\hM$ then $\hT(\partial_r,\partial_r)$ does not depend on $r$ and is a non-constant solution of $(\ref{one})$. Moreover $2\hT$ is the Hessian of the function $A$ defined by $A(r,m):=r^2\hT_{(r,m)}(\partial_r,\partial_r)$.
\end{prop}

In the proof  of proposition \ref{GaMa}, we will need the following two statements; in these statements $X$, $Y$, $Z$   will denote arbitrary vector fields on $M$. We will also 
 denote by the same letters $X,Y,Z$  the lift of these vector fields to $\hM$.

 \begin{fact}[for example, \cite{Leist,Ga,Ma}] \label{fact}
The Levi-Civita connection of $\hg$ is given by 
\begin{equation}\hD_{\wX}\wY={D_XY}-rg(X,Y)\partial_r,\quad \hD_{\partial_r}\partial_r=0, \quad \hD_{\partial_r}\wX=\hD_{\wX}\partial_r=\tfrac{1}{r}\wX.\label{connection}\end{equation}
\end{fact}
{\bf Proof.}  We take a point $(r, m)\in \hM$. Without loss of generality we can assume that $DX(m)=DY(m)=DZ(m)=0$. 
Using $\hg([\wX,\wY],\partial_r)=0$ and $[\partial_r,\wX]=[\partial_r,\wY]=0$, we have
$$2\hg(\hD_{\wX}\wY,\partial_r)=-\partial_r.\hg(\wX,\wY)=-2rg(X,Y).$$
Similarly we have that $2\hg(\hD_{\wX}\wY,\wZ)=r^2g(D_XY,Z)$. It implies the first assertion. 
The two others can be shown the same way. \qed \hip

In the next corollary we  will tautologically  identify $M$ with  $M_{1}:= \{1\}\times M\subset \hM$ (the point $m\in M$ will be identified with $(1,m)\in M_1$). By definition of the cone metric, $\hg_{|M_1}= g$.

\begin{cor}[\cite{Ma}]\label{retour}
 Let $\hT$ be a symmetric parallel $(0,2)-$tensor on $(\hM,\hg)$. Then, $\alpha:= \hT(\partial_r, \partial_r)$ does not depend on $r$ and can be considered therefore as a function on $M$.  Moreover, for every $m\in M$ and every $X,Y,Z\in T_mM$ we have  
\begin{eqnarray}
   2\hT(\partial_r,X)&=&r D\alpha(X) \label{1} \\
   2\hT(X,Y)&=&r^2\left(2\alpha \, g(X,Y)+DD\alpha(X,Y)\right), \label{2}
  \end{eqnarray}
    Moreover,   $\alpha$ is constant if and only if  $\hT$ is proportional to $\hg$. Moreover,  for every $(1,m)\in M_1$ we have  \begin{equation}
  2DT(X,Y,Z)=  -D\alpha\otimes g (Y,X,Z) -  D\alpha\otimes g(Z,X,Y), \label{3} \end{equation}  
  where  $T$ is the restriction of the tensor $\hT$ to $M_{1}$ \  ($\stackrel{\textrm{tautologically}}{\equiv }M$).  
\end{cor}
{\bf Proof.} 
 Since  $\hT$ is parallel, we have
$$ 0=\hD \hT (\partial_r,\partial_r,\partial_r)=\partial_r.\hT(\partial_r,\partial_r)-2\hT(\hD_{\partial_r}\partial_r,\partial_r)\stackrel{\eqref{connection}}{=} \partial_r.\hT(\partial_r,\partial_r)=0.$$ Thus $\hT(\partial_r,\partial_r)$ is a function on $M$. The first statement of corollary \ref{retour} is proved.

Combining   $\hD \hT=0$ with  Fact \ref{fact},  we have:
\begin{equation}
0=\hD \hT(X,\partial_r,\partial_r)=D\alpha(X)-2\hT(\hD_{X}\partial_r,\partial_r)\stackrel{\eqref{connection}}{=}D\alpha(X)-\tfrac{2}{r} \hT(X,\partial_r).\label{up}
\end{equation}
This shows \eqref{1}.
Similarly, using Fact  \ref{fact} and \eqref{1}, we obtain 
$$\begin{array}{rcl}
0=\hD \hT(X,Y,\partial_r)&=&  X.\hT(Y,\partial_r)-\hT(\hD_{X} Y,\partial_r)-\hT(Y,\hD_{X}\partial_r)\\
   &\stackrel{\eqref{connection},\eqref{1}}{=}&\frac{r}{2}\big(DD\alpha(X,Y)+2g(X,Y)\alpha\big) - \tfrac{1}{r}\hT(Y,X).
\end{array}$$
This shows \eqref{2}. If $\alpha=\const$,  $DD\alpha=0$. Then, \eqref{1}, 
 \eqref{2}, and the definition of $\alpha$   implies  $\hT= \alpha \cdot  \hg$.   If $\hT = \const \cdot g$, then $ \hT(\partial_r, \partial_r)= \const\cdot  \hg( \partial_r, \partial_r)= \const. $
  
 Similarly, using Fact  \ref{fact} and already proved parts of corollary \ref{retour},    we obtain (for every $(1,M)\in M_1\subset \hM$)
  $$\begin{array}{rl}
0=\hD \hT(X,Y,Z)=&  X.\hT(Y,Z)-\hT(\hD_{X} Y,Z)-\hT(Y,\hD_{X}Z)\\
   \stackrel{\eqref{connection}}{=}&X.T(Y,Z)- T(D_XY,Z)- T(Y,D_XZ )+     g(X,Y) \hT( Z, \partial_r)  +    g(X,Z) \hT( Y, \partial_r) \\ \stackrel{\eqref{1}}{=}&  DT(X,Y,Z)+ \tfrac{1}{2}\bigl( D\alpha\otimes g (Y,X,Z) +  D\alpha\otimes g(Z,X,Y)  \bigr)  
\end{array}$$
implying \eqref{3}. 
 \qed 

{\bf Proof of proposition \ref{GaMa}.} 
{
Equation \eqref{one} being tensorial we can suppose without loss of generality that $DX=DY=DZ=0$.
}
We 
  set $\hX := \tfrac{1}{r} X$,  $\hY := \tfrac{1}{r} Y$, $\hZ := \tfrac{1}{r} Z$. By  Fact \ref{fact}, we have:
\begin{equation} \label{onemore}
\hD_{\partial_r}\partial_r=\hD_{\partial_r}\hX= \hD_{\partial_r}\hY=\hD_{\partial_r}\hZ=0,\qquad 
\hD_{\hX}\partial_r=\tfrac{1}{r}\hX\quad \mathrm{and} \quad \hD_{\hX}\hY=-\tfrac{1}{r}g(X,Y)\partial_r.\end{equation}

Let $\alpha$ be a solution of $(\ref{one})$,  and $A =r^2\alpha(m)$.  Our first goal is to show that 
$\hD\hD\hD A =0$. 
We have: 
\begin{equation} \label{matveevnew}
\hD A(\hZ)=rD\alpha(Z) \quad \mathrm{and}\quad \hD A(\partial_r)=2r\alpha.\end{equation}
Then,  
\begin{equation}\hD\hD A (\hZ,\partial_r)=\hD\hD A(\partial_r,\hZ)=\partial_r.(r\,D \alpha(Z))-\hD A(\hD_{\partial_r}\hZ)=D\alpha(Z), \label{deriv}\end{equation}
and similarly
\begin{equation}
\hD\hD A(\partial_r,\partial_r)=2\alpha \label{alpha}\end{equation}
Using that, $DY(m)=0$, we get $Y.D A(Z)=DD A(Y,Z)$ and
\begin{equation}\hD\hD A(\hY,\hZ)=\hY.(rD_Z\alpha)+\tfrac{1}{r}g(Y,Z)\hD A(\partial_r)=DD\alpha(Y,Z)+2g(Y,Z)\alpha.\label{hess}\end{equation}

Now we can prove that $\hD \hD \hD  A=0$, we will 
 first   show  that $\hD \hD \hD  A(\partial_r,.,.)=\hD \hD \hD  A(.,\partial_r,.)=\hD \hD \hD A(.,.,\partial_r)=0$:  
$$\begin{array}{lcl}
\hD \hD \hD  A(\partial_r,\hX,\hY) &= & \partial_r . \left( \hD \hD  A(\hX,\hY) \right)\\
& \stackrel{\eqref{hess} }{=} &\partial_r . \left( DD\alpha(Y,Z)+2g(Y,Z)  \right) = 0,\\
 \hD \hD \hD  A(\partial_r,\hX,\partial_r) & =&
\hD \hD \hD  A(\partial_r,\partial_r,\hX) =  \partial_r . \left( \hD \hD  A(\partial_r, \hX) \right) \\
& \stackrel{\eqref{deriv} }{=}& \partial_r . \left( D\alpha(X) \right) = 0,
\\
\hD \hD \hD  A(\partial_r,\partial_r,\partial_r)& = & \partial_r . \left( \hD \hD  A(\partial_r, \partial_r) \right)\\   &\stackrel{\eqref{alpha}} 
=& 2 \partial_r . (\alpha) = 0, \\
 \hD \hD \hD  A(\hX,\partial_r,\partial_r)  &=&
 \hX . \left( \hD \hD  A(\partial_r, \partial_r) \right) - 2 \hD\hD A\left(\hD_{\hX} \partial_r, \partial_r\right) \\ &\stackrel{\eqref{matveevnew},\eqref{alpha} }{=}&  \hX . \left( D\alpha(X) \right) = 0.\\
{
\hD \hD \hD  A(\hX,\partial_r,\hY)} & =&\hD \hD \hD  A(\hX,\hY,\partial_r)\\&\stackrel{\eqref{deriv},\eqref{onemore}}=& { \hX.\hY.\alpha-\frac{1}{r}\hD\hD A(\hX,\hY)-\frac{1}{r}g(X,Y)\hD\hD A(\partial_r,\partial_r)\stackrel{\eqref{alpha},\eqref{hess}}=0}
\end{array}
$$

 The last thing to check is
\begin{eqnarray}
 \hD^3 A(\wX,\hY,\hZ)&=&X.(\hD\hD A(\hY,\hZ))-\hD\hD(\hD_{\wX}\hY,\hZ)-\hD\hD(\hY,\hD_{\wX}\hZ)\nonumber \\ 
                        &=&X.(DD\alpha(Y,Z)+2g(Y,Z)\alpha)+\hD\hD A\left(\tfrac{1}{r}g(X,Y)\partial_r,\hZ\right)+\hD\hD A\left(\tfrac{1}{r}g(X,Z)\partial_r,\hY\right)\nonumber \\
                         &=&DDD \alpha(X,Y,Z)+2g(Y,Z) D\alpha(X)+g(X,Y)D\alpha(Z)+g(X,Z)D\alpha(Y)=0.\nonumber
\end{eqnarray} Thus, $\hD \hD \hD  A=0.$ 
The proposition is proved  in the ``$\Longrightarrow$" direction. 

Let us now prove the proposition in the  ``$\Longleftarrow$" direction.  We take a point $(1,m)\in M_1$. Covariantly differentiating 
\eqref{2} with the help of $D$ and substituting \eqref{3}, we obtain 
$$
\begin{array}{rl}0=&-DT(X,  Y, Z)+  2g(Y,Z)\, X.\alpha+ DDD\alpha(X,Y,Z) \\\stackrel{\eqref{3}}{=}&    2 (D\alpha  \otimes g)(X,Y,Z)+(D\alpha \otimes g)(Y,X,Z)+(D\alpha \otimes g)(Z,X,Y)+DDD \alpha (X,Y,Z).\end{array}
$$

  Comparing corollary \ref{retour} and (\ref{deriv}), (\ref{alpha}), \eqref{hess}  we see that $2\hT$ is the Hessian of $A:=r^2\alpha$.  \qed

Let us recall that 
a pseudo-Riemannian manifold is said to be \slsf{ decomposable} if it possess a non-trivial  parallel non-degenerate (i.e., the restriction of the metric to it is nondegenerate)
 distribution.

By  \cite{Wu},   if a pseudo-Riemannian manifold is decomposable,   then
 the      manifold can be locally written as the  product $(M_1, g_1) \times (M_2 ,g_2)$ 
  of two pseudo-Riemannian manifolds; the tangent space of $M_1$   naturally embedded  in the tangent space of the product is precisely  the parallel distribution. The  tangent space of $M_2$ 
   is the orthogonal complement to the
 parallel distribution, which is itself also a nondegenerate parallel distribution.

Contrarily to the Riemannian case, the existence of a parallel symmetric $(0,2)-$tensor on a pseudo-Rieman\-nian manifold does not imply that the manifold is decomposable. 
It is a consequence of the fact that the self-adjoint endomorphism associated to such a tensor and the metric can not always be simultaneously diagonalized.  However, the situation is more simple for cones  over closed manifolds  as shows the following
\begin{prop}\label{para->decomp}
 Let $(M,g)$ be a closed pseudo-Riemannian manifold and $c=1$. If the equation $(\ref{one})$  has a non-constant solution then $(\hM,\hg)$ is decomposable.
\end{prop}
{\bf Proof.} 
Let $\alpha$ be a non-constant solution of $(\ref{one})$ on $M$. As $M$ is closed there exists two critical points $m_-$ and $m_+$ of $\alpha$  associated to distinct critical values (i.e., $D\alpha(m_{\pm})=0$ and $\alpha(m_-)\neq\alpha(m_+)$).
As $D \alpha(m_{\pm})=0$, it follows from  (\ref{deriv}) and (\ref{alpha}) that, for any $r>0$, $\hD\hD A_{(r,m_{\pm})} (\partial_r,.)$ vanishes on $TM$ and takes the value $2\alpha(m_{\pm})$ on $\partial_r$. It means that
\begin{equation}\hD\hD A_{(r,m_{\pm})} (\partial_r,.)=2\alpha(m_{\pm})g_{(r,m_{\pm})}(\partial_r,.).\label{decadix}\end{equation}
Since the eigenvalues of the parallel tensor are constants, $2\alpha(m_-)\ne2\alpha(m_+)$ are two different eigenvalues of the { field of} self-adjoint endomorphisms  associated to  $\hD\hD A$ at every points.
Since  the  field of self-adjoint endomorphisms  associated to  $\hD\hD A$  is also parallel, 
  its  
  characteristic spaces (=generalized eigenspaces)  provide a parallel orthogonal decomposition    of $T\widehat M$. Indeed, they are clearly nondegenerate;  as there are at least two distinct eigenvalues each characteristic space is non trivial. Then, $(\hM, \hg)$ is decomposable.  \qed  
\hip
Proposition \ref{para->decomp} does not say that a cone over a closed manifold  with interesting holonomy is automatically decomposable. For example, the cone may admit anti-symmetric parallel $(0,2)-$tensors. The reader can consult Alekseevsky et al \cite{Leist} for a more systematic study of the holonomy of cones.

\section{Decomposable cones over closed manifolds and the proof of Theorem \ref{pr}.}\label{decompo}
The goal of this section is to prove 

\begin{prop} \label{deco}   Let  $(M,g)$ be a  closed connected pseudo-Riemannian manifold such that  the  cone $(\hM,\hg)$  is decomposable. Then,   $\hg$ is the Riemannian flat metric, and $g$ is the Riemannian metric of constant curvature $1$. 
\end{prop}

{\bf Proof. }  Let $V_1$ and $V_2:= V_1^{\perp}$ be the  complementary nondegenerate  parallel distribution on $\hM $.  Let $\hT_1$ and $\hT_2$ be the symmetric  $(0,2)$-tensors on $\hM$ defined for $i\in\{1,2\}$ by
$$\hT_i(v,u)=\hg(v_i,u),$$ 
where the $v_i$'s are the factors of the decomposition of $v$ according to the splitting $T\hM=V_1\oplus V_2$.  Clearly, $\hT_1+ \hT_2= \hg$. Since the distributions $V_i$ are parallel,  then the tensors $\hT_i$ are also parallel.

We set $\alpha_i:=\hT_i(\partial_r,\partial_r)$. Since   $\hg(\partial_r,\partial_r)=1$,  we have $\alpha_1+\alpha_2= \hT_1(\partial_r,\partial_r) + \hT_2(\partial_r,\partial_r) = \hg(\partial_r,\partial_r) = 1$.
As in section \ref{section},  we define  the functions $A_1$ and $A_2$ on $\hM$ by $A_i(r,m)=r^2\alpha_i(m)$.

Applying  proposition \ref{GaMa}  to the tensors $\hT_i$ we obtain   
that the following statements hold for every $i=1,2$: 
\begin{itemize}
\item[($\ast$)] $\alpha_i$ is a {non constant} function on $M$ and it is a solution of the equation $(\ref{one})$.
\item[($\ast\ast$)] The Hessian of $A_i$ is $2\hT_i$.
\end{itemize}

Let us prove that 
{\it the only possible  critical values of $\alpha_i$ are $0$ and $1$. 
Moreover, for any $m\in M$, we have $0\leq \alpha_i(m)\leq 1$, the extremal values being reached. }

Since $\alpha_1 + \alpha_2=1$, it is sufficient to prove this statement for $\alpha_1$. 
Let $m\in M$ be  a critical point of $\alpha_1$. As we already saw, at (\ref{decadix}), it implies that  at the point $(r,m)$
$$2\alpha_1(m)g(\partial_r,.)=\hD\hD A_1(\partial_r,.)\stackrel{(\ast\ast)}{=} 2 \hT_1 (\partial_r,.)  .$$ 
Then,  $\partial_r(r,m)$ is an eigenvector  of the self-adjoint endomorphism associated to $\hT_1$, 
and $\alpha_1(m)$ is the eigenvalue of this endomorphism. Since the only  eigenspaces of $\hT_1$ are   $V_1$ (with eigenvalue $1$)  and  $V_2$ (with eigenvalue $0$), then   $\alpha_1(m)=0$ or $\alpha_1(m)=1$. Thus, the only critical values of $\alpha_1$ are $0$ and $1$.  

Since  $M$ is closed, there exists $m_1,m_0\in M^2$  such that $\alpha_1(m_1)=\max_{m\in M} \alpha_1(m)$ and  $\alpha_1(m_0)=\min_{m\in M} \alpha_1(m)$. Then,  $d\alpha_1(m_{1})= d\alpha_1(m_{0})=0$, implying  $\alpha_1 (m_1)=1$ and $\alpha_1(m_0)=0$.

Let us prove that the tensor $\hT_1$ is  nonnegatively definite.  We take the point $(r, m_0)\in \hM$. 
This point is a minimum of  the function $A_1$. Indeed, $A_1(r,m)= r^2 \alpha_1\ge 0$, and $ A_1(r, m_0)= r^2 \alpha_1(m_0)= 0$.  Since $(r,m_0)$ is a minimum, at this point $\hD\hD A_1$ is nonnegatively defined.  Since   $2\hT_1\stackrel{(\ast\ast)}{=} \hD\hD A_1$,  $\hT_1 $ is  nonnegatively defined at the point $(r, m_0)$. Since $\hT_1$ is parallel, it is also  nonnegatively defined at every point of $\hM$. 

Similarly, one can prove that $\hT_2$ is nonnegatively defined: instead of the point $(r, m_0)$ one should take the point $(r, m_1)$ where the function $A_2$ accepts its minimum.

Since $\hg= \hT_1+\hT_2$, it is  also nonnegatively defined. Since it is nondegenerate, it is positively defined, i.e., is a Riemannian metric.  As we recalled in the introduction, the Riemannian version of theorem \ref{pr}  was proved by Gallot \cite{Ga} and Tanno \cite{Ta}. Thus, by Gallot-Tanno Theorem, $g$  has constant curvature equal to $1$, and $\hg$  is the Riemannian flat metric. \qed  

\begin{rema}
The hypothesis of compactness in Proposition \ref{deco} is only use to obtain that the function $\alpha_i$ defined during the proof has a minimum and a maximum. Hence, we could replace the hypothesis of compactness by this weaker one.
\end{rema}

{\bf Proof of theorem \ref{pr}. }  The case $c=0$ was done in section \ref{sec0}. 
 By remark \ref{tttt},  we can assume $c=1$.   By proposition \ref{para->decomp}, the existence of a nonconstant solution of \eqref{one} implies that the cone  $(\hM, \hg)$  is decomposable. By proposition \ref{deco}, $g$  is a Riemannian metric of constant curvature $1$. \qed 

\section{Application I:  geodesic rigidity of Einstein manifolds and projective Obata conjecture.}\label{lichne}
The set of metrics geodesically equivalent (the definition is in \S \ref{his})  to a  metric $g$ is  in one-to-one correspondence with the set of   nondegenerate symmetric $(0,2)-$tensors  $T$ such that for any  vector fields $X,Y,Z$ on  $M$ 
 \begin{equation}\label{basic1} 
 DT(X,Y,Z)= \tfrac{1}{2} \left(D\tr(T)\otimes g (Y,X,Z) +  D\tr(T)\otimes g(Z,X,Y) \right),
  \end{equation}
 where the trace and the covariant derivative are taken according to $g$, see for example \cite[\S 2.2]{einstein} for details (in the ``tensor" notations, the equation \eqref{basic1} reads 
$ T_{ij, k } =  \tfrac{1}{2} \left( T^p_{p,i}  g_{jk} +T^p_{p,j} g_{ik}\right).$)

Since this equation is linear, the space of its solutions is a linear vector space. Its dimension is  called  the \slsf{degree of mobility} of $g$.  

Locally,  the degree of mobility of $g$ coincides with the dimension of the set (equipped with natural topology) of metrics geodesically equivalent to $g$. 

 It is easy  to see that if $\alpha$ is a  solution of the equation $(\ref{one})$ then the tensor defined by \eqref{hess} is a solution of \eqref{basic1}.  Indeed, one can check it directly, or one can use that for   the covariantly-constant by proposition \ref{GaMa}  tensor   $\hT= \hD\hD A$ the tensor \eqref{hess} is precisely the tensor $T$ from corollary \ref{retour}. Then, it satisfies the equation \eqref{3}, which is equivalent to \eqref{basic1}.

 In some cases the reciprocal is true, hence theorem \ref{pr} has the following corollaries.
\begin{cor}  \label{ein} 
Let $g$ be an Einstein (i.e.,  the  Ricci tensor is proportional to $g$) pseudo-Riemannian metric on an $(n>2)-$dimensional  closed connected manifold. Assume that $\bar g$ is geodesically equivalent to $g$, but is not affinely equivalent to  $g$. Then for a certain constant $c\ne 0$ the metric $c\cdot g$ is the Riemannian metric of constant curvature $1$.  
\end{cor}

{\bf Proof.} By  \cite[Corollary 3]{einstein},  if the metric $g$ is Einstein and  if there exists a  geodesically equivalent, but not affine equivalent  metric $\bar g$,   then the equation $(\ref{one})$ admits a non-constant solution. The corollary therefore follows from  theorem \ref{pr}. \qed 

\begin{cor} \label{corobata}  
     Let $g$ be a  pseudo-Riemannian metric on an $(n>1)-$dimensional  closed connected manifold. 
     Then,  if the metric $\bar g$ on $M $  is geodesically equivalent to $g$, but not affinely equivalent to $g$, then the degree of mobility  of $g$ is precisely $2$ or  for certain constants $c\ne 0\ne \bar c $ the metrics 
     $c\cdot g$ and $\bar c\cdot \bar g$ are Riemannian metrics of constant curvature $1$.
\end{cor} 
{\bf Proof.} Assume first that $n=\dim(M)\ge 3$. Under this assumption,  by  \cite{KM}   if the degree of mobility of $g$ is $\ge 3$, then  for every  solution $T$ of \eqref{basic1}, the function $f:= \tr(T)$  is a solution of $(\ref{one})$. More precisely,  \cite[Lemma 3 and Corollary 4]{KM} implies   that in a  neighborhood of almost every point there exists a constant $c$ such that $f$  is a solution of  $(\ref{one})$.  Now, by   \cite[Lemma 7 in \S  2.3.4]{KM}  the constant $c$ is actually universal (implying that the equation  $(\ref{one})$  is fulfilled on the whole manifold). {In this case the result follows therefore from Theorem \ref{pr}.}

Now, the case $n=\dim(M)=2$ follows  from  \cite[Theorem 5.1]{Kio}  for the signatures (+,+) and  (--,--), and from \cite[Corollary 1]{global} for the signature (+, --). 
\qed

Corollary \ref{corobata}  is related to the following classical conjecture: 
\hop
{\bf Projective Obata conjecture.} \emph{
Let $G$ be  a connected Lie group acting on a {closed} connected pseudo-Riemannian or Riemannian manifold $(M,g)$ of dimension $n>1$ by projective transformations. Then it acts by affine transformations or there exists a constant $c\neq0$ such that $(M,c\cdot g)$ is the quotient of a Riemannian round sphere.}
\hop

By corollary \ref{corobata},   we have:
\begin{cor}\label{step}
If $(M^n ,g)$ is a  counter-example to the projective Obata conjecture, then  $n:= \dim(M^n)\ge 3$ and 
the degree of mobility of $g$ is precisely $2$. 
\end{cor}

{\bf Proof.} The existence of a projective nonaffine transformation  for $g$ 
implies the existence of a  metric  that is  geodesically equivalent to $g$, but is not  affine equivalent  to $ g$.    By corollary \ref{corobata},  if  $n=\dim(M^n)\ge 3$,   the degree of mobility of $g$ is  $ 2$.

Now, by  \cite[Theorem 6]{global},  the projective Obata conjecture is true in dimension two.  
  \qed

In the Riemannian case, projective Obata conjecture  was proved in \cite[Theorem 1]{obata} for dimension $2$ and in \cite[Corollary 1]{Ma2} for dimensions $\ge 3$.
The natural idea to prove the  conjecture in the pseudo-Riemannian case is to mimic the Riemannian proof for pseudo-Riemannian metrics.  The (Riemannian) proof contains two  parts:
 \begin{itemize} 
\item[(i)]  proof for the metrics with the degree of mobility $ 2$ (\cite[Theorem 15]{Ma2}), 
\item[(ii)] proof  for the metrics with the degree of mobility $\ge 3$ (\cite[Theorem 16]{Ma2}). 
 \end{itemize} 
We expect that it is possible, though nontrivial, to generalize  (i)  for the  pseudo-Riemannian case. On the other side, one can not expect to generalize (ii) for pseudo-Riemannian metrics, because (ii) is based on Riemannian results that are no more true in a pseudo-Riemannian setting. Hence corollary \ref{step} proves  the part that was, a priori, the most difficult  part of the projective  Obata conjecture  for  pseudo-Riemannian metrics.

 Moreover, the next corollary shows  that 
 the group of projective transformations of a closed  manifold 
coincides with the group of affine transformations, or the group of isometries has codimension one in the 
group of projective transformations.
  \begin{cor} \label{simple} 
Let $(M,g)$  be a closed connected $(n>1)-$dimensional pseudo-Riemannian manifold. Assume that  for no  constant $c\in \mathbb{R}\setminus \{0\}$ the metric $c\cdot g$ is the  Riemannian  metric of  constant  curvature $1$. Then, every projective vector field is an affine vector field, or certain nontrivial linear combination of every two projective vector fields is a Killing vector field.
\end{cor} 
{\bf Proof.}  Indeed,  it is well known (see, for example \cite{Ma2}, or more classical sources acknowledged therein) that a vector field $X$ is projective if the tensor
  \begin{equation} \label{a-2} T:= \wL_X g - \tfrac{1}{n+1} \tr({
 \wL_Xg}) \cdot g \end{equation}  is a solution of \eqref{basic1}, where $\wL_X$ is the Lie derivative with respect to $X$.  Moreover, the projective vector field is affine, if and only if the trace of $T$ is constant.
  
  Suppose the degree of mobility of $g$ is not $2$. Then, corollary  \ref{corobata} implies  that all projective vector fields are actually affine, which is one of the possibilities in corollary \ref{simple}.

  Now, suppose the degree of mobility of $g$ is  precisely $2$. Let $X$ and $Y$ be projective vector fields.  
   We consider  the solutions $T= \wL_X g - \tfrac{1}{n+1} \tr({
\wL_Xg}) \, g $  and $T'= \wL_{Y} g - \tfrac{1}{n+1} \tr({ 
 \wL_Yg}) \cdot g $   of \eqref{basic1}. 
    Since    the degree of mobility is   $2$,   $T$, $T'$, and   $g$  are linearly dependent, i.e., for certain constants $k, k', l$ 
  we have $k\, T +  k'\,T' = l\, g.$
      Since the mapping $$X\mapsto \wL_X g - \tfrac{1}{n+1} \tr(
 {\wL_Xg} )\, g $$ is linear,
      we have 
      $$      \wL_{k\, X + k'\,Y} g - \tfrac{1}{n+1} \tr(
\wL_{k\,X + k' \,Y}g)\,g = l\, g \ \ 
\textrm{implying} \ \  \wL_{k\,X +  k' \,Y} g  = (n+1)l \,g, $$
i.e., $k\, X + k'\,Y$ is a homothety vector field (if $l\ne 0$) or a Killing vector field (if $l=0$). 
      
      Since $M$ is closed, it admits no homotheties implying $k\,X + k'\,Y$  is a Killing vector field. \qed

\section{Application II: Holonomy groups of closed 
constant curvature manifolds.} \label{appendix} 
 Let $X$ be a manifold and $G$ be a Lie group acting analytically on $X$. A  $(G,X)-$\slsf{structure} on a manifold  $M$ is given by an atlas $(U_i,\varphi_i)$ such that each $\varphi_i$ takes values in $X$ and each transition function $\varphi_i\circ\varphi_j^{-1}$ is the restriction of the action of an element of $G$ on $X$.

If $M$ has a $(G,X)-$structure then there exists (see for example \cite[pp. 140,141]{thurston}  for details) a local diffeomorphism $\delta :\widetilde M\rightarrow \widetilde{X}$, where  $\widetilde M$ and $\widetilde{X}$ are the universal covers of $M$ and $X$ respectively,  and a  morphism 
$\rho:\pi_1(M)\rightarrow \widetilde{G}$ (where  $\pi_1(M)$ denotes the fundamental group of $M$,  which acts as the group of deck transformations of  the covering $\widetilde M\rightarrow M$, and $ \widetilde{G}$ denotes the covering of $G$ that acts on $\widetilde{X}$) such that, for any $\gamma \in \pi_1(M)$ and any $\tilde m\in \widetilde M$, we have $\delta (\gamma.\tilde m)=\rho(\gamma).\delta(\tilde m)$.
The map $\delta$ is called the   \slsf{developing map} and the morphism $\rho$ is called the  \slsf{holonomy morphism}. The image of $\pi_1(M)$ with respect to  $\rho$ is called the holonomy group of the $(G,X)-$manifold $M$, it contains a lot of  informations about the geometry of $M$. 

  We denote by $\R^{p+1,q}$ the space $\R^{p+q+1}$ equipped with  the standard pseudo-Euclidean metric of  signature  $(p+1, q)$, and consider the  pseudo-sphere
 $S^{p,q}=\{x\in \R^{p+1,q}\,|\, \langle x,x\rangle=1\}$ and $O(p+1,q)\ltimes\R^{p+q+1}$ the  isometry group of $\R^{p+1,q}$.  We recall that $S^{p,q}$ is  simply connected if and only if $p \neq 1$.

It is well known that every pseudo-Riemannian manifold of signature $(p,q)$ and  constant curvature equal to $1$ is a  manifold with a $(O(p+1,q), S^{p,q})-$structure and that the flat  pseudo-Riemannian manifolds of signature $(p+1,q)$ are the manifolds having a $(O(p+1,q)\ltimes\R^{p+q+1},\R^{p+1,q})$-structure.

The holonomy groups of those $(G,X)-$manifolds are not the usual pseudo-Riemannian holonomy groups (even if they are closely related). For now on we will only consider holonomy groups of $(G,X)-$structures. 
Proposition \ref{deco}   implies
\begin{cor}\label{ghz}
If $q\neq 0$, the action of the  holonomy group of a closed manifold endowed with a $(O(p+1,q), S^{p,q})-$structure
{ 
 (i.e. endowed with a pseudo-Riemannian metric  with constant curvature equal to $1$)
}
 on $\R^{p+1,q}$ does not preserve any non-degenerate splitting. \weg{of $\R^{p+1,q}$.}
\end{cor}
Corollary \ref{ghz} was known under the additional assumption that the manifold is {complete}, see \cite[Fact 2.3]{ze}. Since by \cite{kling}  any constant curvature Lorentz manifold is complete, corollary \ref{ghz} was  also known for closed  manifolds of Lorentz signature.  If $q=0$, the sphere itself is a counterexample.

{\bf Proof of corollary \ref{ghz}.}  In order to simplify the   notation, we will suppose $p\neq 1$,  i.e.,  that $S^{p,q}$ is simply connected. Anyway,  if $p=1$, $M$ is Lorentzian, therefore complete by \cite{kling} and the corollary follows from \cite{ze}. Moreover, it is easy to  adapt  what follows to the case $p=1$.

{ Let $(M,g)$ be a  closed pseudo-Riemannian 
manifold with constant curvature equal to $1$.} We denote by $\widetilde M$ the universal cover of $M$ and  by $\widehat{\widetilde M}$ the cone over its universal cover.
It  is well-known (and follows from Fact \ref{fact}) that the curvature of the cone metric  $\hg= dr^2 + r^2 g$ is given by 
$$\widehat R(X,Y)Z=R(X,Y)Z-g(Y,Z)X+g(X,Z)Y,$$ where $R$ and $\widehat R$ are the curvatures  of $g$ and of  $\hg$. It implies that $(\hM,\hg)$ is flat. 

 We identify  $\widehat{S^{p,q}}$, the cone over the pseudosphere
${S^{p,q}}$,  with $\{x\in \R^{p+1,q}\,|\, \langle x,x\rangle>0\}$.
Let $\delta :\widetilde M \rightarrow S^{p,q}$ be a developing map { of the induced $(O(p+1,q), S^{p,q})-$structure on} $M$.
The  map
$\widehat\delta:\widehat{\widetilde M} \rightarrow  \widehat{S^{p,q}}\weg{=\{x\in \R^{p+1,q}\,|\, \langle x,x\rangle>0\}}$ 
defined by   $\widehat\delta(r,m):=(r,\delta(m))$ is a developing map of 
{ 
the $(O(p+1,q)\ltimes\R^{p+q+1},\R^{p+1,q})$-structure of the flat manifold $\hM$.
} 
The holonomy morphisms associated to $\widehat \delta$ and $\delta$ are clearly the same. 
{ 
We denote them by $\rho$.
}

{ 
Let $T_0$ be a symmetric  parallel $(0,2)-$tensor on  $\R^{p+1,q}$ (i.e. a  symmetric bilinear form) invariant with respect to the holonomy group $\rho(\pi_1(\widehat M))$.
Let $\widehat {\widetilde T}=\widehat\delta^*T_0$ be the pull back of $T_0$ 
by $\widehat\delta$. Let $\gamma$ be an element of $\pi_1(\widehat M)$ seen as the group of deck transformations of the universal covering.
We have $\gamma^*\widehat{ \widetilde T}=\gamma^*(\widehat \delta^* T_0)=\widehat\delta^*(\rho(\gamma)^*T_0)$, but as we supposed that  $\rho(\gamma)^*T_0=T_0$ it implies that $\widehat{ \widetilde T}$ is invariant by the action of $\pi_1(M)$. It means that $\widehat{\widetilde  T}$ is the pull-back of a parallel tensor $\widehat T$ on $\widehat M$.

By propositions \ref{GaMa}, \ref{para->decomp} and \ref{deco},  the tensor $\widehat T$ is proportional to the metric $\widehat g$. Thus $T_0$ also is proportional to the metric of $\R^{p+1,q}$.} It means that the holonomy group of a closed $(O(p+1,q), S^{p,q})-$manifold does not preserve any symmetric bilinear form on  $\R^{p+1,q}$ which is not proportional to the metric. In particular it does not preserve any non-degenerate splitting of $\R^{p+1,q}$.\qed

Note that a stronger version of corollary \ref{ghz}  exists for \emph{flat} pseudo-Riemannian manifolds.  More precisely,  by Goldman et al    \cite{GH}   the holonomy of a closed affine manifold admitting a parallel volume form (for example pseudo-Riemannian and flat) does not preserve any non trivial subspace.

{ As the following example shows, corollary \ref{ghz} is no more true for degenerate splittings.}
 We identify $\R^{2,2}$ with $M(2,\R)$ the space of order $2$ square matrices endowed with the determinant (seen as a quadratic form). The pseudo-sphere is then identified with $SL(2,\R)$. For any  $v\in \R^2\setminus\{0\}$, we define the set $V_v$  by $V_v=\{M\in M(2,\R)\,; \, M.v=0\}$. They are $2$-dimensional totally degenerate subspaces. If $v$ and $w$ are not colinear, we have $\R^{2,2}=V_v\oplus V_w.$ Furthermore, $SL(2,\R)$ clearly acts isometrically by left multiplication on $M(2,\R)$. This action preserves  $V_v$ and $V_w$. Now, let $\Gamma$ be a cocompact lattice in $PSL(2,\R)$. The manifold  $PSL(2,\R)/\Gamma$ is a closed $3$ dimensional anti de Sitter manifold whose  holonomy lies in $SL(2,\R)$ and therefore preserves certain totally degenerate splittings of $\R^{2,2}$.

 However,  up to the authors knowledge,  these examples {(and some of their deformations cf. \cite{salein})} are the only known examples   of closed pseudo-Riemannian manifolds of  constant curvature whose holonomy  preserves a non trivial degenerate subspace.  Moreover, the main proposition  of \cite{ze} is actually that the holonomy  group of a closed anti de Sitter manifold of dimension greater than $3$ is irreducible. \\[.1cm] 

{\bf Acknowledgement:}  V.M.  thanks  Deutsche Forschungsgemeinschaft (Priority Program
1154 � Global Differential Geometry and research training group   1523 --- Quantum and Gravitational Fields) and FSU Jena for partial financial support, and
Vicente Cortes and Dmitri Alekseevsky for useful discussions.

Vladimir S. Matveev\hip
\begin{tabular}{ll}
Address: &Institute of Mathematics, FSU Jena, 07737 Jena Germany\\
E-mail:&{\tt vladimir.matveev@uni-jena.de}
\end{tabular}
\hop
Pierre Mounoud\hip
\begin{tabular}{ll}
 Address: & Universit\'e Bordeaux 1, Institut de Math\'ematiques de Bordeaux,\\
 &351, cours de la lib\'eration, F-33405 Talence, France\\
E-mail:&{\tt pierre.mounoud@math.u-bordeaux1.fr}
\end{tabular}
\end{document}